\documentclass[aap,MSNbibl,seceqn,citesort,dvips]{arximspdf}

\doi{10.1214/11-AAP770}
\volume{22}
\issue{2}
\pubyear{2012}
\firstpage{522}
\lastpage{540}

\makeatletter

\newcommand{\sfont}{\fontsize{8.36pt}{8.36pt}\selectfont}

\newcommand{\CAP}{\operatorname{cap}}

\newtheorem{theorem}{Theorem}[section]
\newtheorem{cor}[theorem]{Corollary}
\newtheorem{lem}[theorem]{Lemma}
\newtheorem{prop}[theorem]{Proposition}

\newproclaim{rem}[theorem]{Remark}

\makeatother

\begin{document}
\begin{frontmatter}

\title{Cover levels and random interlacements\thanksref{T1}}
\runtitle{Cover levels and random interlacements}

\thankstext{T1}{Supported in part by the Grant ERC-2009-AdG 245728-RWPERCRI.}

\begin{aug}
\author[A]{\fnms{David} \snm{Belius}\corref{}\ead[label=e1]{david.belius@math.ethz.ch}\ead[url,label=u1]{http://www.math.ethz.ch/\textasciitilde dbelius/}}
\runauthor{D. Belius}
\affiliation{ETH Z\"{u}rich}
\address[A]{Departement Mathematik\\
ETH Z\"{u}rich\\
CH-8092 Z\"{u}rich\\
Switzerland\\
\printead{e1}\\
\printead{u1}} 
\end{aug}

\received{\smonth{3} \syear{2010}}
\revised{\smonth{12} \syear{2010}}

%
\begin{abstract}
This note investigates cover levels of finite sets in the random
interlacements model introduced in [\textit{Ann. of Math.} (2)
\textbf{171} (2010) 2039--2087], that is, the
least level such that the set is completely contained in the random
interlacement at that level. It proves that as the cardinality of
a set goes to infinity, the rescaled and recentered cover level tends
in distribution to the Gumbel distribution with cumulative distribution
function $\operatorname{exp}(-\operatorname{exp}(-z))$.
\end{abstract}

%
\begin{keyword}[class=AMS]
\kwd[Primary ]{60G50}
\kwd[; secondary ]{82C41}
\kwd{60D05}.
\end{keyword}
\begin{keyword}
\kwd{Cover level}
\kwd{random interlacements}
\kwd{cover time}
\kwd{uncovered set}
\kwd{Gumbel distributional limit}.
\end{keyword}

\end{frontmatter}

\setcounter{section}{-1}
\section{Introduction}\label{intro}

The random interlacements model was introduced\break in~\cite{Sznitman2007}.
It helps to understand the picture left by a simple random walk in
the discrete torus $(\mathbb{Z}/N\mathbb{Z})^{d},d\ge3$, or the discrete
cylinder $(\mathbb{Z}/N\mathbb{Z})^{d-1}\times\mathbb{Z},d\ge3$,
when the walk is run up to times of a certain scale. The random interlacements
are an increasing family of random sets $\mathcal{I}^{u}\subset
\mathbb{Z}^{d}$,
indexed by a parameter $u\ge0$, and for each $u$ the set $\mathcal{I}^{u}$
is, intuitively speaking, the trace the paths whose label is at most
$u$ from a Poisson cloud of labeled doubly infinite paths in $\mathbb{Z}^{d}$
modulo time-shift. By analogy with the concept of random walk cover
times this note introduces the \textit{cover level} of a set by random
interlacements and proves a fine asymptotic limit result for this
quantity. Since random interlacements model random walk in the discrete
torus and cylinder on certain suitable time scales, our result should
eventually lead to a better understanding of the distributional limits
of cover times in these graphs.

We now briefly recall how $\mathcal{I}^{u}$ is constructed. We denote
by $W$ the space of doubly infinite nearest neighbor paths in $\mathbb{Z}^{d}$
that spend finite time in bounded subsets of $\mathbb{Z}^{d}$. We
also introduce the equivalence relation $\sim$ on $W$
by letting $w\sim v$ if $w$ is a time-shift of $v$, that is, if
there exists an $N\in\mathbb{Z}$ such that $w(n)=v(N+n)$ for all
$n\in\mathbb{Z}$. The space $W^{\star}$ of doubly infinite paths
modulo time-shift is defined by $W^{\star}=W/\sim$. The ``Poisson
cloud'' mentioned above is then the interlacement Poisson point process,
that is, a~Poisson point process $\omega=\sum_{i}\delta
_{(w_{i}^{\star},u_{i})}$
on the space $W^{\star}\times[0,\infty)$ with intensity measure given
by the product measure of a certain $\sigma$-finite measure~$\nu$
on~$W^{\star}$ and Lebesgue measure. If $K\subset\mathbb{Z}^{d}$
is finite, the total mass assigned by $\nu$ to the set of trajectories
modulo time-shift that enter $K$ is the capacity of $K$ [see (\ref
{eqDefinitionCapacity})]. If one normalizes the measure $\nu$ on this
set and then
considers for each trajectory modulo time-shift the representative
from $W$ which enters~$K$ for the first time at time $0$, then
$\nu$ corresponds to picking a~position at time $0$ distributed
according to the normalized equilibrium distribution [see~(\ref
{eqdefinitionofEu})], and conditionally on the position at time $0$ the forward
and backward trajectories are, respectively, distributed as simple random
walk and simple random walk conditioned never to reenter $K$. We
refer to Section 1 of \cite{Sznitman2007} for details. We denote
the law governing $\omega$ by $\mathbb{P}$. The random
interlacement~$\mathcal{I}^{u}$ is then defined as
%
%
\begin{equation}\label{eqRandomInterlacementDef}
\mathcal{I}^{u}=\bigcup_{i\dvtx u_{i}\le u}
\operatorname{range}(w_{i}^{\star
}),\ u\ge0, \mbox{ where }\omega=\sum_{i}\delta_{(w_{i}^{\star
},u_{i})},
\end{equation}
and $\operatorname{range}(w^{\star})$ denotes the set of all vertices
of $\mathbb{Z}^{d}$
visited by the path \mbox{$w^{\star}\in W^{\star}$}. The measure~$\nu$
is constructed so that, intuitively speaking, for a value~$u$ related
to the time up to which the random walk in the torus or cylinder is
run, the trace of the random walk ``looks like'' $\mathcal{I}^{u}$
\cite{Windisch2008,Sznitman2009-RWonDTandRI}. The law of the indicator
function of $\mathcal{I}^{u}$ on $\{0,1\}^{\mathbb{Z}^{d}}$ has
a simple characterization (see \cite{Sznitman2007}, (2.16)): it
is the unique law with the property that
%
%
\begin{equation}\label{eqCapacityFormula}
\mathbb{P}(A\cap\mathcal{I}^{u}=\varnothing)=\exp\bigl(-u\cdot
\operatorname
{cap}(A)\bigr)\qquad\mbox{for all finite }A\subset\mathbb{Z}^{d},
\end{equation}
where $\operatorname{cap}(A)$ denotes the capacity of $A$ [see (\ref
{eqDefinitionCapacity})
for the definition].

The sets $\mathcal{I}^{u}$ are naturally increasing in $u$. Taking
inspiration from the concept of cover time of a finite graph one may
consider the \textit{cover level} of a finite set $A\subset\mathbb{Z}^{d}$,
defined as the least level $u$ such that $A$ is completely contained
in $\mathcal{I}^{u}$,
%
%
\begin{equation}\label{eqCoverLevDef}
M(A)=\inf\{u\ge0\dvtx A\subset\mathcal{I}^{u}\}=\max_{x\in
A}U_{x},
\end{equation}
where $U_{x}$ denotes the cover time of the vertex $x$,
%
%
\begin{equation}\label{eqDefOfUx}
U_{x}=\inf\{u\ge0\dvtx x\in\mathcal{I}^{u}\},\ x\in\mathbb{Z}^{d}.
\end{equation}
The main result of this note is the following theorem.
\begin{theorem}[(Rescaled and recentered cover levels have a
distribution close to Gumbel)]
\label{thmMaintheorem}
For any finite nonempty $A\subset\mathbb{Z}^{d}$
we have
%
%
\begin{equation}\label{eqMaintheoremmainstatement}
\sup_{z\in\mathbb{R}}\biggl|\mathbb{P}\biggl(\frac{M(A)}{g(0)}-{\log}|A|\le
z\biggr)-\exp(-e^{-z})\biggr|\le c|A|^{-c_{1}},
\end{equation}
where $c_{1}>0$ is the constant given in (\ref{eqdefintionofc1})
and $g(\cdot)$ is the $\mathbb{Z}^{d}$ Green's function [see (\ref
{eqGreensfunctiondefinition})].
In particular, $\frac{M(A)}{g(0)}-{\log}|A|$ tends in distribution to
the Gumbel distribution, as $|A|$ tends to infinity.
\end{theorem}

We will now give some comments on the scope of the above theorem.
If~$G_{1},G_{2},\ldots$ denotes a sequence of finite graphs whose cardinality
tends to infinity and $C_{N}$ denotes their cover times (i.e., the
first time simple random walk has visited every vertex of the
graph)\vspace*{1pt}
it is sometimes possible (see, e.g., \cite
{MatthewsCoveringProbsForMCs,DevroyeSbihiRWonHighlySymmGraphs})
to show that $\frac{C_{N}}{c|G_{N}|}-{\log}|G_{N}|$ tends in law to
the Gumbel distribution as $N\rightarrow\infty$. Theorem \ref{thmMaintheorem}
should be seen as a result of this flavor. There are, however, simple
families of graphs for which one can show that $\frac
{C_{N}}{{c|G_{N}|\log}|G_{N}|}\rightarrow1$
in probability but for which the finer distributional limit result
remains out of reach to this day. For example, when $G_{N}=(\mathbb
{Z}/N\mathbb{Z})^{d}$, $d\ge3$,
is the discrete torus, it is known that $\frac{C_{N}}{{g(0)|G_{N}|\log}
|G_{N}|}\rightarrow1$
in probability but\vspace*{1pt} only conjectured that $\frac
{C_{N}}{g(0)|G_{N}|}-{\log}|G_{N}|$
tends to the Gumbel distribution (see \cite{aldous-fillbook}, Chapter
7, Section 2.2, pages 22 and~23). In \cite
{Sznitman2009-UBonDTofDCandRI,Sznitman2009-OnDOMofRWonDCbyRI}
and~\cite{TeixeiraWindischOnTheFrag} results are proved that couple
random interlacements and the trace of random walk in the discrete
cylinder $(\mathbb{Z}/N\mathbb{Z})^{d}\times\mathbb{Z},d\ge2$, and
discrete torus $(\mathbb{Z}/N\mathbb{Z})^{d},d\ge3$, in nearly macroscopic
boxes of side length $N^{1-\varepsilon},\varepsilon>0$. In these
works the couplings are used as a ``transfer mechanism'' to allow
one to study the disconnection time of the cylinder and the so-called
``fragmentation'' of the torus by studying a related problem formulated
completely inside the random interlacements model. We believe that
in a similar way the results in the present note for cover levels
of sets in the random interlacements model will lead to progress in
the study of cover times of sets by random walk in the torus and cylinder
(for more on this see Remark \ref{remEndRemark}).

The second result of this note is a corollary of Theorem \ref{thmMaintheorem}.
For any $1\le l\le d$ and $z\in\mathbb{R}$ we let
$B_{N}^{l}=[0,N-1]^{l}\times\{0\}^{d-l}$
and define a sequence $(\mathcal{N}_{N}^{l,z})_{N\ge1}$
of point measures on $\mathbb{R}^{d}$:
%
%
\begin{equation}\label{eqdefinitionofpointproc}
\mathcal{N}_{N}^{l,z}=\sum_{x\in B_{N}^{l}}\delta_{x/N}1_{\{
U_{x}>g(0)\{
{\log}|B_{N}^{l}|+z\}\}},\ N\ge1.
\end{equation}
In other words, $\mathcal{N}_{N}^{l,z}$ collects the points of
$[0,1]^{l}\times\{0\}^{d-l}$
which under scaling correspond to the sites of $B_{N}^{l}$ not yet
covered by the random interlacements at level $g(0)\{{\log}
|B_{N}^{l}|+z\}$.
\begin{cor}[(Convergence of point process of uncovered points to a~homogeneous
Poisson point process)]
\label{corWeakconvergence}
%
%
\begin{equation}\label{eqweakconvergencecorollstatment}
\mathcal{N}_{N}^{l,z}\mbox{ converges in law to }\mathcal
{N}^{l,z}\mbox
{ as }N\rightarrow\infty,
\end{equation}
where $\mathcal{N}^{l,z}$ is a Poisson point process with intensity
$\exp(-z)\lambda_{l}$ and $\lambda_{l}$ is Lebesgue measure on
$[0,1]^{l}\times\{0\}^{d-l}$.
\end{cor}

Incidentally, it follows from this corollary that the last few sites
of $B_{N}^{l}$ to be covered by the random interlacements are ``far
apart,'' at typical distance of about~$N$. This fact is proved in
Proposition \ref{proLastTwoIndians}.

We now comment on the proofs of Theorem \ref{thmMaintheorem} and
Corollary \ref{corWeakconvergence}. For each $x\in\mathbb{Z}^{d}$
the random variable $U_{x}$ is known to have exponential distribution
with parameter $\frac{1}{g(0)}$ [see (\ref{eqDefOfUx}) and
(\ref{eqprobxininterlacement})]. If the $U_{x}$ were independent,
then standard extreme value theory would tell us that the rescaled
and recentered maxima $\frac{M(A)}{g(0)}-{\log}|A|$ tend in distribution
to the Gumbel distribution as $|A|\rightarrow\infty$ (see~\cite
{EmbrechtsKluppelbergMikosch1997}, Example 3.2.7, page 125). However,
in the random interlacements model the $U_{x}$ are not independent,
in fact, there is a long-range correlation; cf. (\ref{eqlongrangecorrelation}).
There is a theory that gives mixing conditions called~$D$ and $D'$
for stationary sequences (\cite{EmbrechtsKluppelbergMikosch1997},
Section 4.4, page~209)
and also similar conditions for stationary random fields \mbox{\cite
{092060037,MR2274809}}
which are sufficient for the rescaled and recentered maxima to converge
in distribution to the Gumbel distribution. When we prove (\ref
{eqlemmacalcwitha})
of Lemma \ref{lemCalculation} we will prove something similar to
$D'$. However, conditions similar to $D$ are difficult to verify
in our context because of the slow decay of the correlation;
cf.~(\ref{eqlongrangecorrelation}). We, therefore, take a different approach
and prove the convergence in distribution directly.

The key to proving Theorem \ref{thmMaintheorem} is to exploit Lemma
\ref{lemRIIndependence},
which says in a quantitative way that in the random interlacements
model spatial separation implies approximate independence. We do this
in Proposition \ref{proCoveringseperatedsetsindependence} by considering
sets $A$ that are ``well separated,'' that is, that consist of isolated
points that are far apart. This spatial separation allows us to use
Lemma \ref{lemRIIndependence} to show that the points of a well-separated
set are covered approximately independently.

We then consider arbitrary nonempty finite sets $A$ and condition
on the subset left uncovered at level ${g(0)(1-\varepsilon)\log}|A|$
for a value of $\varepsilon$ that satisfies $0<\varepsilon\le12c_{1}$.
In Lemma \ref{lemGoodevent} we show that with high probability
this ``uncovered set'' is well separated and has cardinality concentrated
around its expected cardinality, which equals $|A|^{\varepsilon}$.
When the ``uncovered set'' is well separated, Proposition \ref
{proCoveringseperatedsetsindependence},
mentioned above, implies that after level ${g(0)(1-\varepsilon)\log}|A|$
the points of the uncovered set are covered approximately independently.
This allows us to finish the job in the proof of Theorem \ref{thmMaintheorem}
by showing (in a quantitative way) that, when we restrict to a certain
``good event'' that has probability tending to one as $|A|\rightarrow
\infty$,
the cover level $M(A)$ is the maximum of approximately
$|A|^{\varepsilon}$
random variables, which are ``essentially'' independent and exponentially
distributed. It is then not hard to show that if we rescale and recenter
$M(A)$ appropriately, it is close in distribution to the Gumbel distribution,
just as would be the case if the cover levels of the points of $A$
were truly independent.

As alluded to above, random interlacements model the picture left
by random walk in the discrete torus and the discrete cylinder when
run up to times of suitably chosen scales. In this light, the uncovered
set $A_{\varepsilon}$ (when $A=B_{N}^{d}$) can, in particular, be
thought of as a counterpart of the uncovered set in the discrete torus
discussed in \cite{miller-2009,Brummelhuis1991}. See Remarks \ref
{remvarianceremark} and~\ref{remEndRemark} for more on this.

The proof of Corollary \ref{corWeakconvergence} uses Theorem \ref
{thmMaintheorem}
and Kallenberg's theorem (\cite{MR900810}, Proposition 3.22, page 157)
which allow us to verify the convergence of point processes by checking
some straightforward conditions involving convergence of the intensity
measure and the probability that the point measure does not charge
a set.

\section{Notation and a review of random interlacements}\label{sec1}

Constants denoted by $c$ may change from line to line and within
formulas. Unless otherwise indicated, all constants depend only on
the dimension $d$ of $\mathbb{Z}^{d}$. Further dependence on, for
example, parameters $\alpha,\beta$, is denoted by $c(\alpha,\beta)$.
The norm~$|\cdot|$ on $\mathbb{Z}^{d}$ is taken to be the Euclidean
norm. We denote by $|A|$ the cardinality of the set $A$. The notation
$A\subset\subset B$ means that $A$ is a finite subset of $B$. For
two sets $A,B\subset\mathbb{Z}^{d}$ we denote their mutual distance
$\inf_{x\in A,y\in B}|x-y|$ by $d(A,B)$. We define a
path to be a sequence $x_{i},i\ge0$, of elements in $\mathbb{Z}^{d}$
such that $|x_{i+1}-x_{i}|=1$ for $i\ge0$.\vspace*{1pt}

We denote by $W_{+}$ the space of paths in $\mathbb{Z}^{d},d\ge3$,
going to infinity as time goes to infinity. Furthermore, $
(X_{n})_{n\ge0}$
denotes the canonical coordina\-tes,~$\mathcal{W}_{+}$ the $\sigma$-algebra
on $W_{+}$ generated by these coordinates and $\theta
_{n}\dvtx W_{+}\rightarrow W_{+}$
the canonical shift on $W_{+}$. We let $P_{x}$ be the probability
measure on $(W_{+},\mathcal{W}_{+})$ that turns $(X_{n}
)_{n\ge0}$
into a simple random walk starting at~$x$~(sin\-ce for $d\ge3$ the
simple random walk is transient, its law is supported on~$W_{+}$).
If $q\dvtx\mathbb{Z}^{d}\rightarrow[0,\infty)$ then $P_{q}$ denotes
the measure $\sum_{x\in\mathbb{Z}^{d}}q(x)P_{x}$.

Green's function is given by
%
%
\begin{equation}\label{eqGreensfunctiondefinition}
g(x,y)=\sum_{n\ge0}P_{x}(X_{n}=y)\quad\mbox{and}\quad g(\cdot)=g(\cdot
,0).
\end{equation}
Recall the following standard bounds for Green's function (\cite
{MR1117680}, Theorem~1.5.4, page 31):
%
%
\begin{equation}\label{eqgreensfuncestimates}
c|x|^{2-d}\le g(x)\le c|x|^{2-d}.
\end{equation}
Furthermore, by the invariance principle, if $A$ is a set of diameter
at most~$L$ (i.e., $|x-y|\le L$ for all $x,y\in A$) then
%
%
\begin{equation}\label{eqgreensfuncsumestimate}
\sum_{x\in A}g(x)\le cL^{2},
\end{equation}
where we have used that the left-hand side of
(\ref{eqgreensfuncsumestimate})
is bounded by the expected time spent by a random walk in a ball containing
$A$.

For $U\subset\subset\mathbb{Z}^{d}$, we define the entrance time
$H_{U}=\inf\{ n\ge0\dvtx X_{n}\in U\} $
and the hitting time $\tilde{H}_{U}=\inf\{ n\ge1\dvtx X_{n}\in U
\} $.\vadjust{\goodbreak}
The escape probability (or equilibrium measure) $e_{U}\dvtx\mathbb
{Z}^{d}\rightarrow[0,\infty)$
is given by
%
%
\begin{equation}\label{eqdefinitionofEu}
e_{U}(x)=P_{x}(\tilde{H}_{U}=\infty)1_{U}(x)
\end{equation}
and the capacity of $U$ by
%
%
\begin{equation}\label{eqDefinitionCapacity}
\CAP(U)=\sum_{x\in U} e_{U}(x).
\end{equation}
Moreover, for each $y\in K$, we have the equality
%
%
\begin{equation}\label{eqHittingProbabilityForSetDecomposedByLastExit}
P_{x}(H_{K}<\infty)=\sum_{y\in K}g(x,y)e_{K}(y).
\end{equation}

We now recall some further facts about random interlacements. As
mentioned in the \hyperref[intro]{Introduction}, to construct the
random interlacements
one defines on a probability space $(\Omega,\mathcal{A},\mathbb{P})$
the so-called interlacement Poisson point process $\omega$, which
is a Poisson point process on the space $W^{\star}\times[0,\infty)$
of labeled doubly infinite paths modulo time-shift (see the second
paragraph of the \hyperref[intro]{Introduction} for notation) whose
intensity is given
by the product measure of a certain $\sigma$-finite measure $\nu$
and Lebesgue measure. For a detailed construction of the measure
$\nu$ see \cite{Sznitman2007}, Theorem 1.1. In this note we will
only need the existence on the space $(\Omega,\mathcal{A},\mathbb{P})$
of a family of Poisson point processes $\mu_{K,u}$ and $\mu_{K,u,u'}$
on $(W_{+},\mathcal{W}_{+})$, defined for any $K\subset\subset
\mathbb{Z}^{d}$
and any $0\le u\le u'$. Loosely speaking, these point measures on
$W_{+}$ keep track of those doubly infinite paths modulo time-shifts
in $\omega$ that enter $K$, and have labels at most $u$ and labels
between $u$ and $u'$, respectively (i.e., they assign weight 1
to the paths in $W_{+}$ which the double infinite paths modulo time-shifts
induce after their entrance in $K$). The random interlacement $\mathcal
{I}^{u}$,
already defined in terms of the interlacement Poisson process $\omega$
in (\ref{eqRandomInterlacementDef}), can also be constructed
from the point processes $\mu_{K,u}$ (see \cite
{SidoraviciusSznitman2009}, (1.16)),
%
%
\begin{equation}\label{eqdefinitionofIu}
\mathcal{I}^{u} = \bigcup_{K\subset\subset\mathbb
{Z}^{d}}\bigcup
_{w\in\mu_{K,u}}\operatorname{range}(w),\qquad u\ge0.
\end{equation}
[For any point process $\mu$ we write $x\in\mu$ as a shorthand for
$x$ belonging to $\operatorname{Supp}(\mu)$, the support of $\mu$.] For
the definitions and some properties of~$\mu_{K}$, $\mu_{K,u,u'}$ we
refer to \cite{SidoraviciusSznitman2009}, (1.13)--(1.15). We will
need the following facts:
%
%
\begin{eqnarray}\label{eqindependenceofmus}
&&\qquad\begin{tabular}{p{320pt}}
$\mu_{K,u}$ and $\mu_{K,u,u'}$ are independent Poisson point
processes on
$(W_{+},\mathcal{W}_{+})$ with respective intensities $u\cdot
P_{e_{K}}$ and $(u'-u)\cdot P_{e_{K}}$,
\end{tabular}
\\
%
%
\label{eqsumofmus}
&&\qquad\hspace*{6pt}
\mu_{K,u'}=\mu_{K,u}+\mu_{K,u,u'}\quad\mbox{and}\nonumber\\
&&\qquad\hspace*{8.8pt}
\mu_{K,u}=\sum_{i=0}^{m}\delta_{\theta_{H_{K}}(w_{i})}1_{\{
H_{K}(w_{i})<\infty\}}
\mbox{ for } K\subset K'\subset\subset
\mathbb
{Z}^{d}\mbox{ and }\\
&&\qquad\hspace*{8.8pt}
\mu_{K',u}=\sum_{i=0}^{m}\delta_{w_{i}}.\nonumber
\end{eqnarray}
This last compatibility relation also holds with $\mu_{K,u,u'}$ and
$\mu_{K',u,u'}$ replacing~$\mu_{K,u}$ and $\mu_{K',u}$.

From the characterization equation (\ref{eqCapacityFormula}) of the
law of $\mathcal{I}^{u}$ we see that
%
%
\begin{equation}
\label{eqprobxininterlacement}
\mathbb{P}(x\notin\mathcal{I}^{u})  =  \exp\biggl\{ -\frac
{u}{g(0)}\biggr\}
\end{equation}
and
\begin{equation}
\label{eqprobxyininterlacement}
\mathbb{P}(x,y\notin\mathcal{I}^{u}) = \exp\biggl\{ -u\frac
{2}{g(0)+g(x-y)}\biggr\},
\end{equation}
since $\CAP(\{x\})\!=\!\frac{1}{g(0)}$ and $\CAP(\{x,y\})\!=\!\frac{2}{g(0)\!+\!g(x\!-\!y)}$
(see \cite{Sznitman2007}, (1.62) and~(1.64)). As noted in \cite
{Sznitman2007}, (1.68), (\ref{eqprobxyininterlacement}) together
with the bounds on Green's function from (\ref{eqgreensfuncestimates})
implies
%
%
\begin{equation}\label{eqlongrangecorrelation}\qquad
\mathrm{cov}_{\mathbb{P}}\bigl(1_{\{x\in\mathcal{I}^{u}\}},1_{\{y\in
\mathcal
{I}^{u}\}}\bigr)\sim\frac{cu}{|y-x|^{d-2}}\exp(-cu)\qquad\mbox{as
}|x-y|\rightarrow\infty.
\end{equation}
For brevity we write
%
%
\begin{equation}\label{eqdefinitionofuA}
u_{A}(z)=g(0)\{{\log}|A|+z\},
\end{equation}
so that $\{ \frac{M(A)}{g(0)}-{\log}|A|\le z\} =\{
M(A)\le u_{A}(z)\} $.

\section{\texorpdfstring{Proofs of Theorem \protect\ref{thmMaintheorem} and Corollary \protect\ref{corWeakconvergence}}%
{Proofs of Theorem 0.1 and Corollary 0.2}}\label{sec2}

We begin by discussing the overall structure of the proofs of Theorem
\ref{thmMaintheorem}
and Corollary \ref{corWeakconvergence}. The first step is to extract
some independence in the random interlacement model. This is done
in Lemma \ref{lemRIIndependence} which says in a quantitative fashion
that the picture left by the random interlacements in a set $K_{1}$
and the picture left in a set $K_{2}$ are approximately independent
if $K_{1}$ and $K_{2}$ are far apart.

Next we will prove that for well-separated sets, that is, sets consisting
of isolated points that are far apart, the cover levels of the individual
points are approximately independent. This is done in Proposition \ref
{proCoveringseperatedsetsindependence},
which follows easily from Lemma \ref{lemRIIndependence}. If we
were only interested in well-separated sets then Proposition \ref
{proCoveringseperatedsetsindependence}
would be enough to prove convergence to the Gumbel distribution of
the rescaled and recentered cover levels. Intuitively speaking, this
is because the cover level of a~well-separated set is the maximum
of a~set of essentially independent random variables (namely, the cover
levels of the individual points).

But of course we are not dealing only with well-separated sets, but
with arbitrary finite nonempty sets $A$. Therefore, we introduce
the random set~$A_{\varepsilon}$, defined in (\ref{eqdefinitionofAepsilon}),
which consists of all the points of $A$ left uncovered at the level
$g(0){(1-\varepsilon)\log}|A|$ for a parameter $\varepsilon$ such
that
%
%
\begin{equation}\label{eqrangeforepsilon}
0<\varepsilon\le12c_{1},
\end{equation}
where $c_{1}$ is the constant defined in (\ref{eqdefintionofc1}).
For a fixed $\varepsilon\in(0,12c_{1}]$ our methods yield (\ref{eqThisHolds})
(uniformly for all $z$) of which (\ref{eqMaintheoremmainstatement})
is the special case $\varepsilon=12c_{1}$. Whenever $\varepsilon$
appears below, it is always understood to satisfy
(\ref{eqrangeforepsilon}).\vadjust{\goodbreak}

Next we need to show that $A_{\varepsilon}$ is ``well behaved.''
This is done in Lem\-ma~\ref{lemGoodevent} which states that with
probability tending to one, $A_{\varepsilon}$ is well separated and
has cardinality close to $\mathbb{E}|A|$ [which equals
$|A|^{\varepsilon}$
by (\ref{eqExpectedsizeofAepsilon})], or in other words, with
probability tending to one, $A_{\varepsilon}$ belongs to the collection
$G_{A,\varepsilon}$ of ``good sets'' defined in (\ref{eqDefinitionOfG}).

Finally, we turn to the proof of Theorem \ref{thmMaintheorem} in
which the goal essentially is to show that $\mathbb{P}(M(A)\le u_{A}(z))$
is close to the cumulative distribution function of the Gumbel distribution,
that is, close to $\exp(-\exp(-z))$. The major\vspace*{-1pt} step will be to condition
on the set of sites of $A$ not yet\vspace*{2pt} covered at level $(1-\varepsilon
)u_{A}(0)\stackrel{\mbox{\sfont{(\ref
{eqdefinitionofuA})}}}{=}g(0){(1-\varepsilon)\log}|A|$,
that is, on $A_{\varepsilon}$. This will be useful because by (\ref
{eqindependenceofmus})
we have $\mathbb{P}(M(A)\le u_{A}(z)|A_{\varepsilon}=K)=\mathbb
{P}(M(K)\le u_{A}(z)-(1-\varepsilon)u_{A}(0))$.
Next we will show that this latter probability\vadjust{\goodbreak} is close to $\exp(-\exp(-z))$.
Then, multiplying by $\mathbb{P}(A_{\varepsilon}=K)$, summing over
all $K\in G_{A,\varepsilon}$ and using that $A_{\varepsilon}\in
G_{A,\varepsilon}$
with probability tending to one will allow us to show that $\mathbb
{P}(M(A)\le u_{A}(z))$
is close to $\exp(-\exp(-z))$.

The key to proving that $\mathbb{P}(M(K)\le u_{A}(z)-(1-\varepsilon)u_{A}(0))$
is close to $\exp(-\exp(-z))$ is to use the fact that all $K\in
G_{A,\varepsilon}$
are well separated and that $|K|$ is close to $|A|^{\varepsilon}$.
The well-separatedness of $K$ allows us to use Proposition~\ref
{proCoveringseperatedsetsindependence}
to prove that the points of $K$ are covered approximately independently.
Thus $M(K)$ is the maximum of $|K|$ approximately i.i.d. random variables,
meaning that with the correct rescaling and recentering $M(K)$ is
approximately a Gumbel random variable, or in other words, $\mathbb
{P}(M(K)\le u_{K}(z))$
is close to $\exp(-\exp(-z))$. But since $K$ is close to
$|A|^{\varepsilon}$,
it follows that $u_{K}(z)\stackrel{\mbox{\sfont{(\ref
{eqdefinitionofuA})}}}{=}g(0)\{{\log}|K|+z\}$
is close to $u_{A}(z)-(1-\varepsilon)u_{A}(0)\stackrel{\mbox{\sfont
{(\ref{eqdefinitionofuA})}}}{=}g(0)\{{\varepsilon\log}|A|+z\}$.
This will allow us to conclude that $\mathbb{P}(M(K)\le
u_{A}(z)-(1-\varepsilon)u_{A}(0))$
is close to $\mathbb{P}(M(K)\le u_{K}(z))$ and thus also close to
$\exp(-\exp(-z))$.

Corollary \ref{corWeakconvergence} then follows easily from Theorem \ref
{thmMaintheorem}
using Kallenberg's theorem (\cite{MR900810}, Proposition 3.22, page 157).

We begin by stating and proving Lemma \ref{lemRIIndependence}.
The proof is just the calculation leading up to \cite{Sznitman2007}, (2.15),
but it is included here for completeness.

\begin{lem}[(Approximate independence of distant sets in random interlacements)]
\label{lemRIIndependence} Assume $u\ge0$. Let $K_{1},K_{2}\subset
\subset\mathbb{Z}^{d}$
be disjoint sets and let~$B_{1}$, $B_{2}$ be events depending only on
$\mathcal{I}^{u}\cap K_{1}$ and $\mathcal{I}^{u}\cap K_{2}$,
respectively. Then
%
%
\begin{equation}\label{eqRIindependencestatement}
|\mathbb{P}(B_{1}\cap B_{2})-\mathbb{P}(B_{1})\mathbb{P}(B_{2})|\le
cu\frac{\CAP(K_{1})\CAP(K_{2})}{d(K_{1},K_{2})^{d-2}}.
\end{equation}
\end{lem}
\begin{pf}
Decompose $\mu_{K_{1}\cup K_{2},u}=\sum_{n\ge0}\delta_{w_{n}}$ as
follows:
\[
\mu_{K_{1}\cup K_{2},u}=\mu_{1,1}+\mu_{1,2}+\mu_{2,1}+\mu_{2,2},
\]
where
\begin{eqnarray*}
\mu_{1,1} & = & \sum_{n\ge0}\delta_{w_{n}}1_{\{X_{0}\in
K_{1},H_{K_{2}}=\infty\}}, \qquad\mu_{1,2} = \sum_{n\ge0}\delta
_{w_{n}}1_{\{X_{0}\in K_{1},H_{K_{2}}<\infty\},}\\[-3pt]
\mu_{2,1} & = & \sum_{n\ge0}\delta_{w_{n}}1_{\{X_{0}\in
K_{2},H_{K_{1}}<\infty\}}, \qquad\mu_{2,2} = \sum_{n\ge0}\delta
_{w_{n}}1_{\{X_{0}\in K_{2},H_{K_{1}}=\infty\}}.
\end{eqnarray*}
The $\mu_{i,j}$ are simply the restriction of the Poisson point process
$\mu_{K_{1}\cup K_{2},u}$ to disjoint sets and are thus independent
Poisson point processes with respective intensity measures
\begin{eqnarray*}
&&u1_{\{X_{0}\in K_{1},H_{K_{2}}=\infty\}}P_{e_{K_{1}\cup K_{2}}}, \qquad
u1_{\{X_{0}\in K_{1},H_{K_{2}}<\infty\}}P_{e_{K_{1}\cup K_{2}}},\\[-2pt]
&&u1_{\{X_{0}\in K_{2},H_{K_{1}}<\infty\}}P_{e_{K_{1}\cup K_{2}}},
\qquad
u1_{\{X_{0}\in K_{2},H_{K_{1}}=\infty\}}P_{e_{K_{1}\cup
K_{2}}}.
\end{eqnarray*}
There exist measurable functions of point measures $F_{1}$ and $F_{2}$
such that
\[
F_{1}(\mu_{K_{1},u})=1_{B_{1}}\qquad\mbox{a.s.}\quad\mbox{and}\quad
F_{2}(\mu_{K_{2},u})=1_{B_{2}}\qquad\mbox{a.s.}
\]
and thus $\mathbb{P}(B_{1}\cap B_{2})=\mathbb{E}[F_{1}(\mu
_{K_{1},u})F_{2}(\mu_{K_{2},u})]$.
Note that $\mu_{K_{1},u}-\mu_{1,1}-\mu_{1,2}$ is a point process
determined by $\mu_{2,1}$ and thus independent from $\mu_{1,1},\mu
_{1,2},\mu_{2,2}$
and similarly $\mu_{K_{2},u}-\mu_{2,2}-\mu_{2,1}$ is a point process
independent from $\mu_{2,2},\mu_{2,1},\allowbreak\mu_{1,1}$. So we can define
auxiliary point processes $\mu'_{2,1}$ and $\mu'_{1,2}$ such that
$\mu'_{2,1}$ has the same distribution as $\mu_{K_{1},u}-\mu
_{1,1}-\mu_{1,2}$
and $\mu'_{1,2}$ has the same distribution as $\mu_{K_{2},u}-\mu
_{2,2}-\mu_{2,1}$
and $\mu'_{2,1},\mu'_{1,2},\mu_{i,j},1\le i,j\le2$ are\vspace*{-1pt} independent,
so that $\mu_{1,1}+\mu_{1,2}+\mu'_{2,1}\stackrel{\mathrm
{law}}{=}\mu_{K_{1},u}$
and $\mu_{2,2}+\mu_{2,1}+\mu'_{1,2}\stackrel{\mathrm{law}}{=}\mu
_{K_{2},u}$.
Then $\mathbb{P}(B_{1})\mathbb{P}(B_{2})=\mathbb{E}[F_{1}(\mu
_{1,1}+\mu_{1,2}+\mu'_{2,1})F_{2}(\mu_{2,2}+\mu_{2,1}+\mu'_{1,2})]$.
So $|\mathbb{P}(B_{1}\cap B_{2})-\mathbb{P}(B_{1})\mathbb
{P}(B_{2})|$
is bounded above by
%
%
\begin{eqnarray}\label{eqboudbytotalmass}
&&\mathbb{P}(\mu'_{1,2}\ne0\mbox{ or }\mu'_{2,1}\ne0\mbox{ or }\mu
_{1,2}\ne0\mbox{ or }\mu_{2,1}\ne0) \nonumber\\[-9.5pt]\\[-9.5pt]
&&\qquad \le
2\bigl(\mathbb{P}(\mu_{1,2}\ne0)+\mathbb{P}(\mu_{2,1}\ne0)\bigr).\nonumber
\end{eqnarray}
We can bound the probabilities in (\ref{eqboudbytotalmass}) by the
total mass of the intensity measures of the point processes $\mu_{1,2}$
and $\mu_{2,1}$ so that
\begin{eqnarray*}
&&|\mathbb{P}(B_{1}\cap B_{2})-\mathbb{P}(B_{1})\mathbb{P}(B_{2})| \\[-2pt]
&&\qquad\le
2u\bigl(P_{e_{K_{1}\cup K_{2}}}(X_{0}\in K_{1},H_{K_{2}}<\infty)
+P_{e_{K_{1}\cup K_{2}}}(X_{0}\in K_{2},H_{K_{1}}<\infty)\bigr).
\end{eqnarray*}
But note,
\begin{eqnarray*}
P_{e_{K_{1}\cup K_{2}}}(X_{0}\in K_{1},H_{K_{2}}<\infty) & \stackrel
{\mbox{\sfont{(\ref{eqdefinitionofEu})}}}{\le} & \sum_{x\in
K_{1}}e_{K_{1}}(x)P_{x}(H_{K_{2}}<\infty)\\[-2pt]
& \stackrel{\mbox{\sfont{(\ref
{eqHittingProbabilityForSetDecomposedByLastExit})}}}{=} & \sum_{x\in
K_{1},y\in K_{2}}e_{K_{1}}(x)g(x,y)e_{K_{2}}(y)\\[-2pt]
& \stackrel{\mbox{\sfont{(\ref{eqgreensfuncestimates}), (\ref
{eqDefinitionCapacity})}}}{\le} & cd(K_{1},K_{2})^{2-d}\CAP
(K_{1})\CAP(K_{2}).
\end{eqnarray*}
Applying a similar calculation for $P_{e_{K_{1}\cup K_{2}}}(X_{0}\in
K_{2},H_{K_{1}}<\infty)$
we get~(\ref{eqRIindependencestatement}).
\end{pf}\eject

We are now ready to prove Proposition \ref{proCoveringseperatedsetsindependence}
which says that the points of well-separated sets are covered approximately
independently. When we use it in the proof of Theorem \ref{thmMaintheorem}
we will use a value for the parameter $\lambda$ which depends on
$\varepsilon$.
\begin{prop}[(The points of well-separated sets are covered almost
independently)]
\label{proCoveringseperatedsetsindependence}Let $\lambda>0$
be a parameter and let $A\subset\subset\mathbb{Z}^{d}$ be nonempty
and such that $|x-y|\ge|A|^{({2+\lambda})/({d-2})}$ for all distinct
$x,y\in A$. Then for $u\ge0$ we have
%
%
\begin{equation}\label{eqcoveringserperatedsetsindependence}
\bigl|\mathbb{P}\bigl(M(A)\le u\bigr)-[\mathbb{P}(U_{0}\le u)]^{|A|}\bigr|\le
cu|A|^{-\lambda}.
\end{equation}
\end{prop}
\begin{pf}
Fix an arbitrary $x\in A$ and let $B_{1}=\{ M(A\setminus\{x\})\le
u\} $
and $B_{2}=\{ M(\{x\})\le u\} =\{ U_{x}\le u\} $.
Applying Lemma \ref{lemRIIndependence} we get
\begin{eqnarray*}
&&\bigl|\mathbb{P}\bigl(M(A)\le u\bigr)-\mathbb{P}(U_{0}\le u)\mathbb
{P}\bigl(M(A\setminus\{x\})\le u\bigr)\bigr|\\
&&\qquad \le cu\frac{\CAP(A\setminus\{
x\})\CAP(\{x\})}{d(\{x\},A\setminus\{x\})^{d-2}}
\stackrel{\mbox{\sfont{(\ref{eqDefinitionCapacity}), (\ref
{eqgreensfuncestimates})}}}{\le} cu\frac{|A|}{|A|^{2+\lambda
}}\\
&&\qquad=cu|A|^{-1-\lambda}.
\end{eqnarray*}
Now applying the same step $|A|-1$ more times, with the other elements
of $A$ substituted for $x$, and the appropriate subsets of $A$
substituted for $A\setminus\{x\}$, and using the triangle inequality
we get
\[
\bigl|\mathbb{P}\bigl(M(A)\le u\bigr)-[\mathbb{P}(U_{0}\le u)]^{|A|}\bigr|\le
cu|A||A|^{-1-\lambda}=cu|A|^{-\lambda}.
\]
\upqed\end{pf}
\begin{rem}
Assume $A_{1},A_{2},\ldots$ is a sequence of sets with
$|A_{i}|\rightarrow
\infty$
as $i\rightarrow\infty$ that are well separated in the sense that
they satisfy the hypothesis of Proposition \ref
{proCoveringseperatedsetsindependence}
for some fixed $\lambda>0$. Then convergence in distribution of the
rescaled and recentered cover levels of the $A_{i}$ to the Gumbel
distribution follows, since, in the notation of (\ref{eqdefinitionofuA}),
$[\mathbb{P}(U_{0}\le u_{A_{i}}(z))]^{|A_{i}|}=(1-\frac
{\exp(-z)}{|A_{i}|})^{|A_{i}|}$
tends to $\exp(-\exp(-z))$ and the right-hand side of (\ref
{eqcoveringserperatedsetsindependence})
tends to zero as $|A_{i}|\rightarrow\infty$. This observation is
a first step on the way to arbitrary sets.
\end{rem}

We now define the constant $c_{1}$ that appears in Theorem \ref
{thmMaintheorem}:
%
%
\begin{equation}
\label{eqdefintionofc1}
c_{1} = \frac{1}{4}\min\biggl(\frac{1}{14}\frac{d-2}{d-1},\frac
{c_{2}}{9-c_{2}}\biggr)
\end{equation}
and
\begin{equation}
\label{eqdefintionofc2}
c_{2} = P_{0}(\tilde{H}_{0}=\infty).
\end{equation}
Since the random walk is transient in $\mathbb{Z}^{d}$ for $d\ge3$,
we have $c_{2}>0$, so that $c_{1}>0$.
\begin{rem}
Since $P_{0}(\tilde{H}_{0}<\infty)\sim\frac{1}{2d}$ as
$d\rightarrow
\infty$
we see that $\frac{c_{2}}{9-c_{2}}\rightarrow\frac{1}{8}>\frac{1}{14}$
as $d\rightarrow\infty$, so for large $d$ we have $c_{1}=\frac
{1}{56}\frac{d-2}{d-1}$.
Actually it can be shown that $\frac{c_{2}}{9-c_{2}}>\frac{1}{14}$
for all $d\ge3$, and hence, $c_{1}=\frac{1}{56}\frac{d-2}{d-1}$ for
all $d\ge3$. We do\vspace*{1pt} not include the details, but to do this, one uses
the expression for $g(0)$ (when $d=3$) in terms of an integral given
in \cite{MR0088110}, (4.1), and the\vspace*{2pt} explicit computation of this
integral (scaled by a factor $\frac{1}{3}$) from \cite{MR0001257},
together\vspace*{-1pt} with the trivial bound $K(k)\le\frac{\pi}{2}\frac{1}{\sqrt
{1-k^{2}}}$
on $K(k)$, the complete elliptic integral of the first kind.
\end{rem}

The next result, Lemma \ref{lemCalculation}, encapsulates a calculation
used in Lemma~\ref{lemGoodevent} to prove that with probability
tending to one $A_{\varepsilon}$, defined in (\ref{eqdefinitionofAepsilon}),
is a~``good set,'' that is, belongs to the collection
$G_{A,\varepsilon}$
from (\ref{eqDefinitionOfG}). We recall that we tacitly assume
$0<\varepsilon\le12c_{1}$.
\begin{lem}
\label{lemCalculation}For any $A\subset\subset\mathbb{Z}^{d}$ let
%
%
\begin{equation}\label{eqdefinitionofAepsilon}
A_{\varepsilon}=\{x\in A\dvtx U_{x}>g(0){(1-\varepsilon)\log}|A|\}
\end{equation}
denote the set of points of $A$ not yet covered at level
$(1-\varepsilon
)u_{A}(0)\stackrel{\mbox{\sfont{(\ref
{eqdefinitionofuA})}}}{=}g(0){(1-\varepsilon)\log}|A|$.
Then for all $A\subset\subset\mathbb{Z}^{d}$ and $b\ge1$ we have
%
%
\begin{equation}
\label{eqlemmacalcwithouta}
\sum\mathbb{P}(x,y\in A_{\varepsilon}) \le |A|^{2\varepsilon
}+c(\varepsilon)|A|^{-{\varepsilon}/{3}}
\end{equation}
and
%
%
\begin{equation}
\label{eqlemmacalcwitha}
\sum\mathbb{P}(x,y\in A_{\varepsilon}) \le c(\varepsilon
)b^{d}|A|^{-{\varepsilon}/{3}},
\end{equation}
where\vspace*{1pt} the first sum is over all distinct $x,y\in A$ and the second
sum is over all $x,y\in A$ such that $0<|x-y|<b|A|^{{1}/({2(d-1)})}$.
\end{lem}
\begin{pf}
We will prove that
%
%
\begin{equation}\label{eqlemmalcalcgeneral}\qquad
\sum_{x,y\in A,0<|x-y|<a}\mathbb{P}(x,y\in A_{\varepsilon})\le\min
(|A|^{2\varepsilon},c|A|^{2\varepsilon-1}a^{d})+c(\varepsilon
)|A|^{-{\varepsilon}/{3}}.
\end{equation}
This\vspace*{-2pt} implies (\ref{eqlemmacalcwitha}) by taking $a=b|A|^{{1}/({2(d-1)})}$
and noting that $2\varepsilon-1+\frac{d}{2(d-1)}\le-\frac
{1}{3}\varepsilon$
because $1-\frac{d}{2(d-1)}=\frac{1}{2}\frac{d-2}{d-1}\stackrel
{\mbox
{\sfont{(\ref{eqdefintionofc1})}}}{\ge}28c_{1}\stackrel{\mbox
{\sfont
{(\ref{eqrangeforepsilon})}}}{\ge}\frac{7}{3}\varepsilon$.
Also\vspace*{1pt} (\ref{eqlemmacalcwithouta}) follows from (\ref{eqlemmalcalcgeneral})
by letting $a\rightarrow\infty$. We begin by splitting the sum in
(\ref{eqlemmalcalcgeneral}) into
\[
I_{1} = \sum_{x,y\in A,0<|x-y|\le({\log}|A|)^{2}}\mathbb
{P}(x,y\in
A_{\varepsilon})
\]
and
\[
I_{2} = \sum_{x,y\in A,({\log}|A|)^{2}<|x-y|<a}\mathbb{P}(x,y\in
A_{\varepsilon}).
\]
To bound $I_{1}$ we note
%
%
\begin{eqnarray}\label{eqboundonI1}
I_{1} & \stackrel{\mbox{\sfont{(\ref
{eqprobxyininterlacement})}}}{=} &
\sum_{x,y\in A,0<|x-y|\le({\log}|A|)^{2}}\exp\biggl(-{(1-\varepsilon)\log}
|A|\frac
{2g(0)}{g(0)+g(x-y)}\biggr)\nonumber\hspace*{-30pt}\\[-9.5pt]\\[-9.5pt]
& \le& c({\log}|A|)^{2d}|A|^{1-2(1-\varepsilon)/
({1+g(e_{1})/g(0)})},\nonumber
\end{eqnarray}
where in the inequality we have used that $\frac
{g(e_{1})}{g(0)}=\mathbb
{P}_{e_{1}}(H_{0}<\infty)\ge\mathbb{P}_{x}(H_{0}<\infty)=\frac{g(x)}{g(0)}$
for all $x\in\mathbb{Z}^{d},x\ne0$. The exponent of $|A|$ equals
$1-(1-\varepsilon)\frac{2}{2-c_{2}}$ since $\frac
{g(e_{1})}{g(0)}\stackrel{\mbox{\sfont{(\ref
{eqdefintionofc2})}}}{=}1-c_{2}$.
The definition of $c_{1}$ and (\ref{eqrangeforepsilon}) immediately
imply that $\varepsilon\le12c_{1}\le3\frac{c_{2}}{9-c_{2}}<3\frac
{c_{2}}{8-c_{2}}$
and rearranging gives $c_{2}>\frac{8\varepsilon}{\varepsilon+3}$.
Plugging this in we have $1-(1-\varepsilon)\frac
{2}{2-c_{2}}<1-(1-\varepsilon)\frac{2}{2-8\varepsilon/(\varepsilon
+3)}=1-(1-\varepsilon)\frac{2\varepsilon+6}{6-6\varepsilon}=1-\frac
{1}{3}(\varepsilon+3)=-\frac{\varepsilon}{3}$.
So $I_{1}$ is bounded above by $c(\varepsilon)|A|^{-{\varepsilon}/{3}}$.
To bound $I_{2}$ first note that using the elementary inequality
$\frac{1}{1+x}\ge1-x$, we find
%
%
\begin{eqnarray}\label{eqexpofgreensfunctionsum}\qquad
I_{2} & \stackrel{\mbox{\sfont{(\ref
{eqprobxyininterlacement})}}}{=} &
\sum_{x,y\in A,({\log}|A|)^{2}<|x-y|<a}\exp\biggl(-2{(1-\varepsilon)\log}
|A|\frac{1}{1+{g(x-y)}/{g(0)}}\biggr)\nonumber\\[-2pt]
& \le& \sum_{x,y\in A,({\log}|A|)^{2}<|x-y|<a}\exp\biggl(-2{(1-\varepsilon
)\log}|A|\biggl(1-\frac{g(x-y)}{g(0)}\biggr)\biggr)\\[-2pt]
& \le& |A|^{-2(1-\varepsilon)}\sum_{x,y\in A,({\log}
|A|)^{2}<|x-y|<a}\exp\bigl({c\log}|A|g(x-y)\bigr).\nonumber
\end{eqnarray}
Now note that $g(x-y)\stackrel{\mbox{\sfont{(\ref
{eqgreensfuncestimates})}}}{\le}c|x-y|^{2-d}\le c({\log}|A|)^{-2}$
for $|x-y|\ge({\log}|A|)^{2}$ so the quantity in the exponential
in (\ref{eqexpofgreensfunctionsum}) is bounded. Thus we can
conclude that (\ref{eqexpofgreensfunctionsum}) itself is bounded
by
\begin{eqnarray*}
&&|A|^{-2(1-\varepsilon)}\sum_{x,y\in A,({\log}|A|)^{2}<|x-y|<a}\bigl(1+{c\log}
|A|g(x-y)\bigr)\\[-2pt]
&&\qquad \le
\min(|A|^{2\varepsilon},c|A|^{2\varepsilon
-1}a^{d})+{c|A|^{-2(1-\varepsilon)}\log}|A|\sum_{y\in A}\sum_{z\in A-y}g(z).
\end{eqnarray*}
We have no control on the diameter of $A$. But let $A-y=\{
a_{1},a_{2},\ldots,a_{n}\}$
with $|a_{1}|\le|a_{2}|\le\cdots\le|a_{n}|$, $|A|=n$, and let
$b_{1},b_{2},\ldots$
be an enumeration of~$\mathbb{Z}^{d}$ with
$|b_{1}|\le|b_{2}|\le\cdots.$
Then by Green's function estimates in (\ref{eqgreensfuncestimates})
\[
g(a_{i})\le c|a_{i}|^{2-d}\le c|b_{i}|^{2-d}\le cg(b_{i}).
\]
So $\sum_{z\in A-y}g(z)\le c\sum_{i=1}^{n}g(b_{i})\le c|A|^{{2/d}}$
since the diameter of $\{b_{1},\ldots,b_{n}\}$ is bounded by $cn^{1/d}$.
Hence,
\begin{eqnarray*}
I_{2} & \le& \min(|A|^{2\varepsilon},c|A|^{2\varepsilon
-1}a^{d})+{c|A|^{-2(1-\varepsilon)}\log}|A||A|^{{2/d}+1}\\
& \le& \min(|A|^{2\varepsilon},c|A|^{2\varepsilon-1}a^{d})+{c\log
}|A||A|^{-{10\varepsilon/9}},
\end{eqnarray*}
where we have used that $-2(1-\varepsilon)+\frac{2}{d}+1=2\varepsilon
-\frac{d-2}{d}\le2\varepsilon-\frac{2}{3}\frac{d-2}{d-1}\stackrel
{\mbox
{\sfont{(\ref{eqdefintionofc1})}}}{\le}2\varepsilon-\frac
{112}{3}c_{1}\stackrel{\mbox{\sfont{(\ref
{eqrangeforepsilon})}}}{\le
}-\frac{10}{9}\varepsilon$.
Combining this with $I_{1}\le c(\varepsilon)|A|^{-\varepsilon/3}$
then gives (\ref{eqlemmalcalcgeneral}).
\end{pf}

Our next task is to use the above lemma to prove that with high probability
$A_{\varepsilon}$ is ``well behaved.''
\begin{lem}[(The good event is likely)]
\label{lemGoodevent}For $A\subset\subset\mathbb{Z}^{d}$ and
$\varepsilon$
as in (\ref{eqrangeforepsilon}) let
%
%
\begin{eqnarray}\label{eqDefinitionOfG}
G_{A,\varepsilon} & = & \bigl\{K\subset A\dvtx\bigl||K|-|A|^{\varepsilon}
\bigr|\le|A|^{2\varepsilon/3},K\ne\varnothing\mbox{ and
}\nonumber\\[-8pt]\\[-8pt]
&&\hspace*{4.3pt} |x-y|\ge(2^{1/\varepsilon}|A|)^{{1}/({2(d-1)})}\mbox{
for all distinct }x,y\in K\bigr\}
\nonumber
\end{eqnarray}
denote the collection of subsets of $A$ that are well separated and
close in cardinality to $|A|^{\varepsilon}$. Then for all $A\subset
\subset\mathbb{Z}^{d}$
one has
%
%
\begin{equation}\label{eqgoodeventstatement}
\mathbb{P}(A_{\varepsilon}\notin G_{A,\varepsilon})\le c(\varepsilon
)|A|^{-{\varepsilon}/{3}}.
\end{equation}
\end{lem}
\begin{pf}
The following two statements together clearly imply (\ref
{eqgoodeventstatement}):
%
%
\begin{eqnarray}
\label{eqgoodeventpart1}\qquad
\mathbb{P}\bigl(\exists x,y\in A_{\varepsilon}\mbox{ s.t.
}0<|x-y|<(2^{1/\varepsilon}|A|)^{{1}/({2(d-1)})}\bigr) & \le&
c(\varepsilon)|A|^{-{\varepsilon}/{3}},\\
\label{eqgoodeventpart2}
\mathbb{P}\bigl(\bigl||A_{\varepsilon}|-|A|^{\varepsilon}
\bigr|>|A|^{{2\varepsilon}/{3}}\bigr) & \le& c(\varepsilon)|A|^{-
{\varepsilon}/{3}}.
\end{eqnarray}
By the union bound (\ref{eqgoodeventpart1}) follows directly
from (\ref{eqlemmacalcwitha}) with $b=\break2^{({1}/{\varepsilon})
({1}/({2(d-1)}))}\ge1$.
To prove (\ref{eqgoodeventpart2}) note that
%
%
\begin{equation}\label{eqExpectedsizeofAepsilon}
\mathbb{E}|A_{\varepsilon}|=|A|\mathbb{P}\bigl(U_{0}>(1-\varepsilon
){g(0)\log}|A|\bigr)\stackrel{\mbox{\sfont{(\ref
{eqprobxininterlacement})}}}{=}|A|^{\varepsilon}.
\end{equation}
So by Chebyshev's inequality,
%
%
\begin{equation}\label{eqchebyshev}
\mathbb{P}\bigl(\bigl||A_{\varepsilon}|-|A|^{\varepsilon}
\bigr|>|A|^{{2\varepsilon}/{3}}\bigr) \le\frac{\mathbb
{E}|A_{\varepsilon}|^{2}-|A|^{2\varepsilon}}{|A|^{
{4\varepsilon}/{3}
}}.
\end{equation}
But
\begin{eqnarray*}
\mathbb{E}|A_{\varepsilon}|^{2} & = & \sum_{x,y\in A}\mathbb
{P}(x,y\in
A_{\varepsilon})\\
& = & \sum_{x\in A}\mathbb{P}(x\in A_{\varepsilon})+\sum_{x,y\in
A,x\ne y}\mathbb{P}(x,y\in A_{\varepsilon})\\
& \stackrel{\mbox{\sfont{(\ref{eqlemmacalcwithouta})}}}{\le} &
|A|^{\varepsilon}+|A|^{2\varepsilon}+c(\varepsilon)|A|^{-\varepsilon/3}.
\end{eqnarray*}
Plugging this bound for $\mathbb{E}|A_{\varepsilon}|^{2}$ into (\ref
{eqchebyshev})
then gives (\ref{eqgoodeventpart2}).
\end{pf}
\begin{rem}
\label{remvarianceremark}In bounding the numerator of the right-hand
side of (\ref{eqchebyshev}) we showed that the variance of
$|A_{\varepsilon}|$
is bounded\vspace*{-1pt} from above by $|A|^{\varepsilon}+c(\varepsilon
)|A|^{-\varepsilon/3}$.
The\vadjust{\goodbreak} inequality\vspace*{1pt} $\mathbb{E}|A_{\varepsilon}|^{2}=|A|^{\varepsilon
}+\sum
_{x\ne y}\mathbb{P}(x,y\in A_{\varepsilon})\stackrel{\mbox{\sfont
{(\ref{eqprobxyininterlacement})}}}{\ge}|A|^{\varepsilon
}+|A|(|A|-1)|A|^{-2(1-\varepsilon)}$
gives a matching lower bound and proves that\break $\operatorname
{Var}|A_{\varepsilon}|\sim|A|^{\varepsilon}$
as $|A|\rightarrow\infty$. In particular, if $A=B_{N}^{d}$, then
$\operatorname{Var}|A_{\varepsilon}|\sim N^{d\varepsilon}$.
As might be expected given the connection between random interlacements
and random walk in the discrete torus (see \cite
{Windisch2008,TeixeiraWindischOnTheFrag}),
this agrees with the value which was found in the theoretical physics
paper \cite{Brummelhuis1991} for the variance of the number of points
of the torus $(\mathbb{Z}/N\mathbb{Z})^{d},d\ge3$, not covered by
random walk run up to time $(1-\varepsilon)g(0){N^{d}\log} N^{d}$ (see
\cite{Brummelhuis1991}, (3.15), (3.9)). Note that the time
$(1-\varepsilon)g(0){N^{d}\log} N^{d}$
is a fraction $1-\varepsilon$ of the ``typical'' cover time
$g(0){N^{d}\log} N^{d}$
(see \cite{aldous-fillbook}, Chapter 7, Section 2.2, page 22,
Corollary 24)
of the torus, just as $(1-\varepsilon){g(0)\log} N^{d}$ is a fraction
$1-\varepsilon$ of the ``typical'' cover level ${g(0)\log} N^{d}$
of $B_{N}^{d}$ by random interlacements.
\end{rem}

We are now ready to prove the main theorem. Recall the definition
of~$u_{A}(z)$ from (\ref{eqdefinitionofuA}).
\begin{pf*}{Proof of Theorem \ref{thmMaintheorem}}
If $z\le-\frac{1}{4}{\varepsilon\log}|A|$ and $|A|>1$, then
\[
\mathbb{P}\bigl(M(A)\le u_{A}(z)\bigr)\stackrel{\mbox{\sfont{(\ref
{eqdefinitionofAepsilon})}}}{\le}\mathbb{P}(A_{\varepsilon
/4}=\varnothing
)\stackrel{\mbox{\sfont{(\ref{eqDefinitionOfG})}}}{\le}\mathbb
{P}(A_{\varepsilon/4}\notin G_{A,\varepsilon/4})\stackrel{\mbox
{\sfont{(\ref{eqgoodeventstatement})}}}{\le}c(\varepsilon
)|A|^{-\varepsilon/12}.
\]
Also, $\exp(-\exp(-z))\le\exp(-|A|^{\varepsilon/4})\le
|A|^{-\varepsilon/4}$
so that
%
%
\begin{equation}\label{eqThisHolds}
\bigl|\mathbb{P}\bigl(M(A)\le u_{A}(z)\bigr)-\exp(-\exp(-z))\bigr|\le c(\varepsilon
)|A|^{-\varepsilon/12}.
\end{equation}
Furthermore, if $z\ge{\log}|A|$ then
\begin{eqnarray*}
\mathbb{P}\bigl(M(A)>u_{A}(z)\bigr) & \le& \mathbb{P}\bigl(A\not\subset\mathcal
{I}^{{2g(0)\log}|A|}\bigr)\\
& \stackrel{\mbox{\sfont{(\ref{eqprobxininterlacement})}}}{\le} &
|A|\exp(-{2\log}|A|)=|A|^{-1}.
\end{eqnarray*}
Also, for such $z$ we have $\exp(-\exp(-z))\ge\exp(-|A|^{-1})\ge1-|A|^{-1}$
so (\ref{eqThisHolds}) also holds for $z\ge{\log}|A|$. It thus
remains to show (\ref{eqThisHolds}) for $z\in(-\frac
{1}{4}{\varepsilon
\log}|A|$, ${\log}|A|)$,
and in what follows we assume $z$ to be in this range.

Let $\mu_{1}=\mu_{A,(1-\varepsilon)u_{A}(0)},\mu_{2}=\mu
_{A,(1-\varepsilon)u_{A}(0),u_{A}(z)}$
and $\mu_{3}=\mu_{A,u_{A}(z)}$. If we fix any $K\in G_{A,\varepsilon}$
then $\{A_{\varepsilon}=K\}$ is simply the event $E_{1}=\{A\setminus
K=A\cap\bigcup_{w\in\mu_{1}}\operatorname{range}(w)\}$
[recall from the remark above (\ref{eqdefinitionofIu}) that $w\in u_{1}$
means\break $w\in\operatorname{Supp}(\mu_{1})]$. Furthermore, $\{M(A)\le
u_{A}(z)\}$
is simply the event $\{A\subset\break\bigcup_{w\in\mu_{3}}\operatorname
{range}(w)\}$,
and since by (\ref{eqsumofmus}) we have $\mu_{1}+\mu_{2}=\mu_{3}$,
the intersection $E_{1}\cap\{M(A)\le u_{A}(z)\}$ coincides with
$E_{1}\cap\{A\subset(A\setminus K)\cup\bigcup_{w\in\mu
_{2}}\operatorname{range}(w)\}=E_{1}\cap\{K\subset\bigcup_{w\in\mu
_{2}}\operatorname{range}(w)\}$.
Denote the last event in the latter intersection by $E_{2}$, and
recall that by (\ref{eqindependenceofmus}) we have that $\mu_{1}$
and $\mu_{2}$ are independent, so that $\mathbb{P}(E_{1}\cap
E_{2})=\mathbb{P}(E_{1})\times\mathbb{P}(E_{2})$.
Also by (\ref{eqindependenceofmus}) the point process $\mu
_{A,u_{A}(z)-(1-\varepsilon)u_{A}(0)}\stackrel{\mbox{\sfont{(\ref
{eqdefinitionofuA})}}}{=}\mu_{A,g(0)\{{\varepsilon\log}|A|+z\}}$
has the same law as $\mu_{2}$, so
\[
\mathbb{P}(E_{2})=\mathbb
{P}\biggl(K\subset
\bigcup_{w\in u_{A,g(0)\{{\varepsilon\log}|A|+z\}}}\operatorname
{range}(w)\biggr)=\mathbb{P}\bigl(M(K)\le g(0)\{{\varepsilon\log}|A|+z\}\bigr).
\]
It follows that (\ref{eqmarkovpropofRI}) holds for all $K\in
G_{A,\varepsilon}$
by noting that both the right- and the left-hand side equal $\mathbb
{P}(E_{1}\cap E_{2})/\mathbb{P}(E_{1})$.
%
%
\begin{equation}\label{eqmarkovpropofRI}\qquad
\mathbb{P}\bigl(M(A)\le u_{A}(z)|A_{\varepsilon}=K\bigr) = \mathbb
{P}\bigl(M(K)\le
g(0)\{{\varepsilon\log}|A|+z\}\bigr).
\end{equation}
Then consider some $K\in G_{A,\varepsilon}$. Define $\lambda=\frac
{1}{2\varepsilon}\frac{d-2}{d-1}-2\stackrel{\mbox{\sfont{(\ref
{eqrangeforepsilon}), (\ref{eqdefintionofc1})}}}{\ge}\frac
{14}{6}-2=\frac{1}{3}$
and note that $\frac{\lambda+2}{d-2}\varepsilon=\frac{1}{2(d-1)}$,
so that for distinct $x,y$ in $K$ we have
\begin{eqnarray*}
|x-y|\stackrel{\mbox{\sfont{(\ref{eqDefinitionOfG})}}}{\ge
}(2^{1/\varepsilon}|A|)^{{1}/({2(d-1)})}=(2|A|^{\varepsilon
})^{
({2+\lambda})/({d-2})} & \stackrel{\mbox{\sfont{(\ref
{eqDefinitionOfG})}}}{\ge} & |K|^{({2+\lambda})/({d-2})}.
\end{eqnarray*}
We can now use Proposition \ref{proCoveringseperatedsetsindependence}
to get that for $|A|>1$,
%
%
\begin{eqnarray}\label{eqTheorem1B}
&&\bigl|\mathbb{P}\bigl(M(K)\le g(0)\{{\varepsilon\log}|A|+z\}\bigr)
-\mathbb{P}\bigl(U_{0}\le
g(0)\{{\varepsilon\log}|A|+z\}\bigr)^{|K|}\bigr|\nonumber\\[-2pt]
&&\qquad
\stackrel{z\le{\log}|A|}{\le} {c\log}|A||K|^{-\lambda}\\[-2pt]
&&\hspace*{11.7pt}\qquad \le {c(\varepsilon)\log}|A||A|^{-{\varepsilon}/{3}},\nonumber
\end{eqnarray}
since $|K|^{-\lambda}\stackrel{\mbox{\sfont{(\ref
{eqDefinitionOfG})}}}{\le}(|A|^{\varepsilon}-|A|^{2\varepsilon
/3})^{-\lambda}\le c(\varepsilon)|A|^{-\lambda\varepsilon}\le
c(\varepsilon)|A|^{-{\varepsilon}/{3}}$.
Furthermore, using again that $||A|^{\varepsilon}-|K||\le|A|^{
{2\varepsilon}/{3}}$
and (\ref{eqprobxininterlacement}) we see that
%
%
\begin{eqnarray}\label{eqTheorem1C}
\biggl(1-\frac{e^{-z}}{|A|^{\varepsilon}}\biggr)^{|A|^{\varepsilon
}+|A|^{2\varepsilon/3}} & \le& \mathbb{P}\bigl(U_{0}\le g(0)\{{\varepsilon
\log}|A|+z\}\bigr)^{|K|}\nonumber\\[-9.5pt]\\[-9.5pt]
& \le& \biggl(1-\frac{e^{-z}}{|A|^{\varepsilon}}
\biggr)^{|A|^{\varepsilon}-|A|^{2\varepsilon/3}}.\nonumber
\end{eqnarray}
But note that if $|A|\ge c(\varepsilon)$ one has the inequality
\[
\biggl|\exp(-e^{-z})-\biggl(1-\frac{e^{-z}}{|A|^{\varepsilon}}
\biggr)^{|A|^{\varepsilon}\pm|A|^{2\varepsilon/3}}\biggr| \stackrel{z\ge
-({\varepsilon}/{4}){\log}|A|}{\le} c|A|^{-\varepsilon/12}.
\]
Combining (\ref{eqmarkovpropofRI}), (\ref{eqTheorem1B})
and (\ref{eqTheorem1C}) with the above formula yields that
if $|A|\ge c(\varepsilon)$,
%
%
\begin{equation}\label{eqCloseToGumbelConditioned}\qquad
\bigl|\mathbb{P}\bigl(M(A)\le u_{A}(z)|A_{\varepsilon}=K\bigr)-\exp(-e^{-z})
\bigr|\le c(\varepsilon)|A|^{-\varepsilon/12}.
\end{equation}
Multiplying by $\mathbb{P}(A_{\varepsilon}=K)$, and summing over
all $K\in G_{A,\varepsilon}$, we see that
\[
\bigl|\mathbb{P}\bigl(M(A)\le u_{A}(z),A_{\varepsilon}\in G_{A,\varepsilon
}\bigr)-\exp(-e^{-z})\mathbb{P}(A_{\varepsilon}\in G_{A,\varepsilon})
\bigr|\le c(\varepsilon)|A|^{-\varepsilon/12}.
\]
Finally, two applications of (\ref{eqgoodeventstatement}) gives
us that (\ref{eqThisHolds}) holds for $|A|\ge c(\varepsilon)$
and $z\in(-\frac{1}{4}{\varepsilon\log}|A|,{\log}|A|)$. Thus (\ref
{eqThisHolds})
holds for all $z\in\mathbb{R}$, and (\ref{eqMaintheoremmainstatement})
follows by taking $\varepsilon=12c_{1}$.\vspace*{-3pt}
\end{pf*}

We now use Theorem \ref{thmMaintheorem} to prove Corollary \ref
{corWeakconvergence}
which states that the point process of uncovered points of a box converges
to a homogeneous Poisson point process.\vspace*{-3pt}
\begin{pf*}{Proof of Corollary \ref{corWeakconvergence}}
We drop the superscripts on $\mathcal{N}_{N}^{l,z}$ and $\mathcal{N}^{l,z}$
to lighten the notation. By Kallenberg's\vadjust{\goodbreak} theorem (see \cite{MR900810},
Proposition 3.22, \mbox{page}~157) it suffices to check that
%
%
\begin{eqnarray}
\label{eqKallenbergcondition1}
\lim_{N\rightarrow\infty}\mathbb{E}\mathcal{N}_{N}(I) & = &
\mathbb
{E}\mathcal{N}(I)\mbox{ for all }I\in\mathcal{J}
\end{eqnarray}
and
\begin{equation}
\label{eqKallenbergcondition2}
\lim_{N\rightarrow\infty}\mathbb{P}\bigl(\mathcal{N}_{N}(I)=0\bigr) =
\mathbb{P}\bigl(\mathcal{N}(I)=0\bigr)\mbox{ for all }I\in\mathcal{J},
\end{equation}
where $\mathcal{J}=\{\mbox{Finite unions of open rectangles }\prod
_{i=1}^{d}(a_{i},b_{i})\mbox{ in }\mathbb{R}^{d}\}$.
To veri-\break fy~(\ref{eqKallenbergcondition1}) simply note that
\begin{eqnarray*}
\mathbb{E}\mathcal{N}_{N}(I) & = & \sum_{x\in NI\cap
B_{N}^{l}}\mathbb
{P}\bigl(U_{x}>g(0)\{{\log}|B_{N}^{l}|+z\}\bigr)\\[-2pt]
& \stackrel{\mbox{\sfont{(\ref{eqprobxininterlacement})}}}{=} &
\frac
{|NI\cap B_{N}^{l}|}{|B_{N}^{l}|}\exp(-z)\rightarrow\lambda
_{l}(I)\exp
(-z)\qquad\mbox{as }N\rightarrow\infty.
\end{eqnarray*}
Since $\mathbb{E}\mathcal{N}(I)=\lambda_{l}(I)\exp(-z)$ by the definition
of $\mathcal{N}$ we have proved (\ref{eqKallenbergcondition1}).
To verify (\ref{eqKallenbergcondition2}) note
\begin{eqnarray*}
\mathbb{P}\bigl(\mathcal{N}_{N}(I)=0\bigr) & = & \mathbb{P}\bigl(M(NI\cap
B_{N}^{l})\le g(0)\{{\log}|B_{N}^{l}|+z\}\bigr)\\[-2pt]
& = & \mathbb{P}\bigl(M(NI\cap B_{N}^{l})\le g(0)\{{\log}|NI\cap
B_{N}^{l}|+z'\}\bigr),
\end{eqnarray*}
where $z'=\log\frac{|B_{N}^{l}|}{|NI\cap B_{N}^{l}|}+z$. So applying
Theorem \ref{thmMaintheorem} we get that for all $N\ge1$,
\[
\biggl|\mathbb{P}\bigl(\mathcal{N}_{N}(I)=0\bigr)-\exp\biggl(-\frac{|NI\cap
B_{N}^{l}|}{|B_{N}^{l}|}\exp(-z)\biggr)\biggr|\le c|NI\cap B_{N}^{l}|^{-c_{1}}.
\]
By taking the limit $N\rightarrow\infty$ and using that $\frac
{|NI\cap
B_{N}^{l}|}{|B_{N}^{l}|}\rightarrow\lambda_{l}(I)$
we get
\[
\mathbb{P}\bigl(\mathcal{N}_{N}(I)=0\bigr)\rightarrow\exp(-\lambda
_{l}(I)\exp(-z)).
\]
But $\mathbb{P}(\mathcal{N}(I)=0)=\exp(-\lambda_{l}(I)\exp(-z))$,
so (\ref{eqKallenbergcondition2}) follows.\vspace*{-3pt}
\end{pf*}

Corollary \ref{corWeakconvergence} has the following interesting implication.
\begin{prop}
\label{proLastTwoIndians}Let the random vector $(X_{1},X_{2},\ldots
,X_{N^{l}})$
be the sites of $B_{N}^{l}$ ordered by the level at which they are
covered (where we use, e.g., the lexicographic order when several
sites are covered at the same level) so that
\[
M(B_{N}^{l})=U_{X_{1}}\ge U_{X_{2}}\ge\cdots\ge U_{X_{N^{l}}}.
\]
Then for all $k\ge2$,
%
%
\begin{equation}\label{eqlasttwoindiansmainstatement}\qquad
\lim_{\delta\rightarrow0}\limsup_{N\rightarrow\infty}\mathbb
{P}(\exists
1\le i<j\le k\mbox{ such that }|X_{i}-X_{j}|\le\delta N)=0,
\end{equation}
or in other words, the last sites of $B_{N}^{l}$ to be covered by
the random interlacements are separated, at typical distance of order
$N$.
\end{prop}
\begin{pf}
Fix a $\delta>0$ and let $f\dvtx\mathbb{R}^{d}\rightarrow[0,\delta^{-1}]$
be a continuous function such that $f(x)=\delta^{-1}$ when\vadjust{\goodbreak} $|x|\le
\delta$
and $f(x)=0$ when $|x|\ge2\delta$. Consider the sum $\sum_{x,y\in
\mathcal{N}_{N},x\ne y}f(x-y)=\mathcal{N}_{N}\otimes\mathcal
{N}_{N}(f(\cdot-\cdot))-f(0)\mathcal{N}_{N}([0,1]^{d})$
[recall that $\mathcal{N}_{N}=\mathcal{N}_{N}^{l,z}$ and $\mathcal
{N}=\mathcal{N}^{l,z}$
depend on $l$ and $z$ and also the remark above~(\ref{eqdefinitionofIu})
about the notation $x\in\mathcal{N}_{N}$]. We have that $\mathcal{N}_{N}$
tends weakly to $\mathcal{N}$, so the product $\mathcal{N}_{N}\otimes
\mathcal{N}_{N}$
tends weakly to $\mathcal{N}\otimes\mathcal{N}$, so that for all
$z\in\mathbb{R}$,
%
%
\begin{eqnarray}\label{eqintermsofPPP}
\lim_{N\rightarrow\infty}\mathbb{E}\biggl[\sum_{x,y\in\mathcal
{N}_{N},x\ne
y}f(x-y)\biggr] & = & \mathbb{E}\bigl[\mathcal{N}\otimes\mathcal{N}\bigl(f(\cdot
-\cdot
)\bigr)-f(0)\mathcal{N}([0,1]^{d})\bigr]\nonumber\hspace*{-30pt}\\[-9.5pt]\\[-9.5pt]
& = & \mathbb{E}\biggl[\sum_{x,y\in\mathcal{N},x\ne y}f(x-y)\biggr].\nonumber
\end{eqnarray}
Let $\tilde{\mathcal{N}}$ be a homogeneous Poisson point process
on $\mathbb{R}^{l}\times\{0\}^{d-l}$ (identified with $\mathbb{R}^{l}$)
with intensity measure $\exp(-z)\lambda_{l}$. Recall that the Palm
measure of $\tilde{\mathcal{N}}$ (viewed as a point process on
$\mathbb{R}^{l}$)
is simply $\exp(-z)$ times the law of $\tilde{\mathcal{N}}+\delta
_{\{0\}}$
(\cite{NeveuProcPonct}, Chapter~2, Exercise 3). So by the definition
of the Palm measure (\cite{NeveuProcPonct}, Chapter 2, Theorem~II.4)
we get
\begin{eqnarray*}
\mathbb{E}\biggl[\sum_{x,y\in\tilde{\mathcal{N}}}1_{\{x\in
[0,1]^{l},y\ne x\}
}f(x-y)\biggr] & = & \exp(-z)\int_{[0,1]^{l}}\mathbb{E}[\mathcal{\tilde
{N}}(f)]\,dx\\[-2pt]
& = & \exp(-2z)\int_{\mathbb{R}^{l}}f(x)\,dx\\[-2pt]
&\le& c(z)\delta^{l-1}.
\end{eqnarray*}
The left-hand side of the above equality is an upper bound for the
right-hand side of (\ref{eqintermsofPPP}). So for any $z\in\mathbb{R}$
we find that
\begin{eqnarray*}
&&\lim_{\delta\rightarrow0}\limsup_{N\rightarrow\infty}\mathbb
{P}\bigl(\exists
x\ne y\mbox{ in }\operatorname{Supp}(\mathcal{N}_{N})\mbox{ such that
}|x-y|\le\delta\bigr)\\[-2pt]
&&\qquad \le
\lim_{\delta\rightarrow0}\limsup_{N\rightarrow\infty}\mathbb
{P}\biggl(\sum
_{x,y\in\mathcal{N}_{N},x\ne y}f(x-y)\ge\delta^{-1}\biggr)\\[-2pt]
&&\qquad \le
\lim_{\delta\rightarrow0}\delta\limsup_{N\rightarrow\infty
}\mathbb
{E}\biggl(\sum_{x,y\in\mathcal{N}_{N},x\ne y}f(x-y)\biggr)  = 0.
\end{eqnarray*}
Finally, we have the inequality
\begin{eqnarray*}
&&\mathbb{P}(\exists1\le i<j\le k\mbox{ such that }|X_{i}-X_{j}|\le
\delta N) \\[-2pt]
&&\qquad \le
\mathbb{P}\bigl(\exists x\ne y\mbox{ in } \operatorname{Supp}(\mathcal
{N}_{N})\mbox{ such that }|x-y|\le\delta\bigr) + \mathbb{P}\bigl(\mathcal
{N}_{N}([0,1]^{d})<k\bigr).
\end{eqnarray*}
So taking, in order, the limits $N\rightarrow\infty,\delta
\rightarrow0$
and $z\rightarrow-\infty$ (and noting that $\lim_{z\rightarrow
-\infty
}\mathbb{P}(\mathcal{N}([0,1]^{d})<k)=0$)
we get (\ref{eqlasttwoindiansmainstatement}).
\end{pf}

We finish with a remark about the possible applicability of our results
to the study of random walk cover times, a comment about the connection
between the uncovered set $A_{\varepsilon}$ and the uncovered set
in the discrete torus and an open question about whether our results
can be generalized.\vadjust{\goodbreak}
\begin{rem}
\label{remEndRemark}(1) Using our results and the known connection
between random interlacements and simple random walk in the discrete
cylinder (see \cite
{Sznitman2009-OnDOMofRWonDCbyRI,Sznitman2009-UBonDTofDCandRI})
it should be possible to determine the finer asymptotic behavior of
the cover time of the cylinder's zero level $(\mathbb{Z}/N\mathbb
{Z})^{d-1}\times\{0\},d\ge3$,
by random walk. Present knowledge states that the cover time is asymptotic
to $N^{2d(1+o(1))}$ (see \cite{Sznitman2006}, Theorem 1).

(2) As already explained in the \hyperref[intro]{Introduction} in the
paragraph after
the statement of Theorem \ref{thmMaintheorem}, it is tempting to use
the coupling result from~\cite{TeixeiraWindischOnTheFrag} together
with Theorem \ref{thmMaintheorem} to devise a proof of the conjecture
that $\frac{C_{N}}{g(0)N^{d}}-\log N^{d}$ tends in law to the Gumbel
distribution, where $C_{N}$ denotes the cover time of the discrete
torus of side length $N$ and dimension $d\ge3$.

(3) In this note the uncovered set $A_{\varepsilon}$ is studied
as one step in proving fine results about the covering of sets by
random interlacements. The corresponding uncovered set in the torus
(cf. Remark \ref{remvarianceremark}) has been studied for its own
sake. Further illustrating the connection between random interlacements
and random walk in the torus, $A_{\varepsilon}$ and the uncovered
set in the torus share some properties. Other than the agreement of
the variance of the cardinality of the uncovered sets mentioned in
Remark~\ref{remvarianceremark}, \cite{miller-2009} also shows that
in the torus the uncovered set is (in a certain sense) well separated
(\cite{miller-2009}, Lemma 6.4), a result similar in spirit to
our Lemma \ref{lemGoodevent}.

(4) Random interlacements can be constructed for any infinite graph
on which simple random walk is transient (see \cite
{TeixeiraIntPercOnTransWeighGraphs}).
It is an open question whether a result like Theorem \ref{thmMaintheorem}
can be proved for random interlacements on more general graphs. It
seems plausible that on a transient graph $G$ such that Green's function
decays ``fast enough'' and such that $a=\operatorname{cap}(\{x\})$ is
independent of $x\in G$ one can use the same method to prove that
for sequences of finite sets $A\subset G$
\[
\frac{M(A)}{a^{-1}}-{\log}|A|\stackrel{\mathrm{law}}{\rightarrow
}\mbox
{Gumbel distribution, as }|A|\rightarrow\infty,
\]
where $M(A)$ is the cover level of $A$.
\end{rem}


\section*{Acknowledgments}

The author would like to thank Alain-Sol Sznitman for suggesting the
problem and supervising his research and Augusto Teixeira for useful
discussions.


%

%
\printaddresses

\end{document}